# Sustainable Ecosystem Planning Based on Discrete Stochastic Dynamic Programming and Evolutionary Game Theory

Yuan Zhang, *Graduate Student Member, IEEE*

*Abstract*—This paper proposed a discrete stochastic dynamic programming (SDP) model for sustainable ecosystem (SE) planning of the Loess Plateau in Northwestern, China, and analyzed the ecological resource planning by the evolutionary game model in the decision-making process. The main objective is to explore a new approach of SE planning from a viewpoint of discrete SDP and evolutionary game theory, with a specific application in the area of ecological resource planning such as water management problems. In contrast to previous work, the proposed SDP method focuses on the transition probability matrix of the ecosystem in a statistic sense, and uses the DP algorithm to obtain the optimal ecological resource planning strategies among multi-subsystems, then analyzes impacts of decision between different users. Firstly, the application background and the concept of SE planning are introduced. Then, a brief overview of existing theory for analyzing sustainable ecosystem is presented. Furthermore, a SDP-based mathematical model and its application to water resource planning of central areas of Loess Plateau are presented as an example. Finally, supplementary analysis of impacts between different users in SE planning as a game playing is provided.

*Index Terms*—The Loess Plateau, sustainable ecosystem planning, water resource planning, discrete stochastic dynamic programming, transition probability matrix, evolutionary game theory.

## I. Introduction

SINCE THE 20TH CENTURY, ecosystems (including biosphere) have been faced with threats under the impacts of climate and humankind together. Sustainability is an important target of developing nature ecosystem. The different patterns of resource utilization could directly influenced ecosystem heath [1], which however, could threatened the existence and development of humankind itself in turn, especially in economically undeveloped and sparsely population areas. A typical example of those threats comes from decades of *high-input* ecological practices (insecticides, logging, overexploitation, etc.) have resulted in increased short-term yields at the price of long-term degradation of ecosystems [2]. Thus, sustainable ecosystem (SE), as a focal point of modern ecology research, faces server challenges. The reasonable adjustment, resource planning, and construction of SE should favorably accorded with the functions of natural, economic and social compound system, and promoted sustainable utilization of resources, such as land, water, vegetation [3]-[5]. Thus, the actual concept of SE planning should be addressed in the framework of systematic and dynamic evolutionary processes.

Fortunately, a large body of literature on control and optimization of dynamic systems has been developed. There have been attempts to apply some of the concepts and techniques to problems in resource ecology [6]. Since then, system analysis and models of ecosystems have progressed from a purely descriptive to an interpretative level [6]. Up to now, numerous models of rangelands, forests, crops soil, antelope conservation and landscape reconstruction have been developed and incorporated with goal-oriented management approaches to ecosystem practices [2]. Specifically, research of SE planning in the past decades has been initiated to find ways to reverse the aforementioned detrimental trend [7], [8]. Edwards emphasized the need to study interactions among the components of integrated ecosystems, recognizing that these systems are complex and are influenced by many factors which are usually difficult to define and quantify. Diverse flows of energy, mass and information into, within, and out of the system mix, and in many cased amalgamate. Moreover, the impact of social and economic circumstances on the system should not be ignored [5]. It is not generally trivial to build a model to meet the needs of SE planning to accommodate these diverse influences, especially at the small-scale regional level or in an integrated economic system [5].

Nevertheless, a fundamental weakness in many earlier SE models is that they use strictly deterministic and quantitative approaches to describe systems that are full of uncertainty and only qualitatively understood. In fact, what is generally needed in developing SE models are methodologies that can represent the multi-agent components of the system and their dynamic interactions in an analytical form using a reasonable number of equations and parameters [2]. To meet these demands, several models based on decision theory framework were proposed and have made important contributions to resource planning of SE through the use of more qualitative techniques like hierarchical ranking methods [9] and quantitative techniques such as stochastic dynamic programming (SDP) [10], [11]. Moreover, SDP has been used to find the best ways of harvesting metapopulations [1], releasing a biological control agent [12], maintaining ecosystem diversity [13], and translocating individuals between two species to ensure persistence [14]. Although the above models considered uncertainties and focused on the analysis of overall resource planning problem among multi-subsystems, they ignore impacts of dynamic relationship among them, namely *evolutionary game relations*. An illustrative example can be expressed as collective deer hunting game of how to work together to achieve sustainable

This work was supported by the Dean's Fellowship Program of Boston University. Y. Zhang is with the Division of System Engineering, College of Engineering, Boston University, 15 Saint Mary's Street, Brookline, MA 02446 USA. E-mail: yzboston@bu.edu.

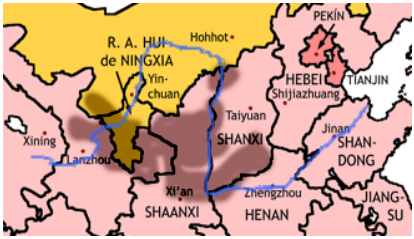

Fig. 1. Graphic view of Loess Plateau, where is shaded. (cited from Wikipedia [15]).

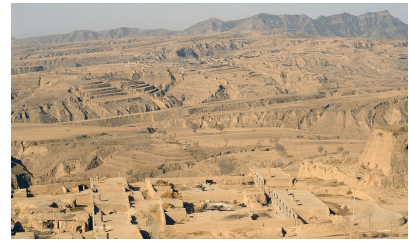

Fig. 3. Graphic view of Loess Plateau near Hunyuan in Shanxi Province (cited from Wikipedia [15]).

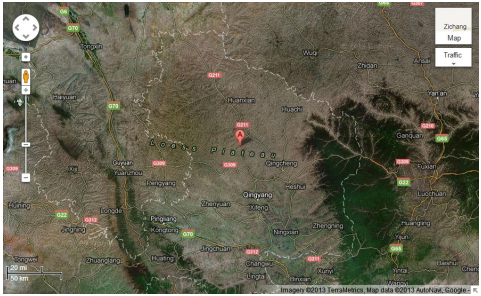

Fig. 2. Graphic view of Loess Plateau from Google Maps (cited from http://maps.google.com).

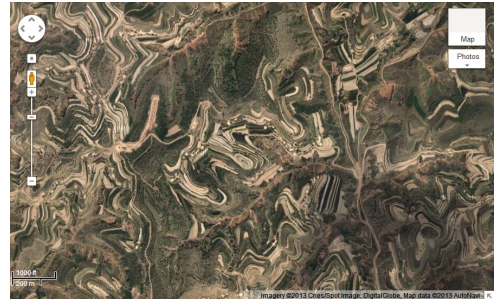

Fig. 4. Topography of Loess Plateau from Google Maps (cited from http://maps.google.com).

use of resources among different users.

In fact, previous works were mainly focus on economically developed and densely populated areas. They neglected regions with adverse weather conditions, such as arid and semi-arid areas with low natural productivity and high sensitivity to disturbance. For example, Loess Plateau in Northwestern China, is a highly erosion-prone soil that is susceptible to the forces of wind and water; in fact, the soil of this region has been called the "most highly erodible soil on earth" [15], which has already threatened the existence and development of SE within this area. Since most parts of the area belong to semi-arid zones and desert distributions, which is of continental climate and influenced very little by the monsoon. The average annual precipitation is about 400mm. Severe shortage and imbalance of water resource along with drought has already become bottleneck of the economic/social development, and is rising to be the primary constraint. Hence it is urgent to build models of sustainable management and utilization of water resource for Loess Plateau by considering both SDP and evolutionary game theory.

In this paper, a stochastic multi-dimensional model of dynamic programming is proposed based on discrete SDP, and analyzes the problem of water resource planning of the Loess Plateau by the evolutionary game model in decision game playing processes. Unlike the previous studies, it focuses on both ecological resource planning problem among multi-subsystems and impacts of their dynamic game relationship. For readability and effectiveness of exposition, we illustrate the method and its salient features using the central region of Loess Plateau (named as *Yulin* area) as an example. The analysis is of help to us in doping out the optimal water resource allocation policies, and more significantly, some overall strategies for sustainable management of Loess Plateau.

The paper is organized as follows. Section II gives brief introduction of Loess Plateau from ecological viewpoints. Section III presents a brief overview of the complexity and game theory. In Section IV, a model based on SDP is formulated for the Loess Plateau ecological resource planning, which uses data of Yulin area so as to analyze optimal water resource planning problems as an example. The results of the empirical analysis provide \valuable references to the SE planning of the Loess Plateau. Section V proposes the evolutionary game model which describes the dynamic decisions of how to cooperate together to achieve sustainable modes of resources among different users. It is a useful supplement to SDP model, and could accordingly lead to some insightful strategic recommendations of SE planning in the Loess Plateau. The paper also contains one Appendix, which address the transition probability matrix.

## II. BRIEF OVERVIEW OF LOESS PLATEAU

The Loess plateau covers an extensive region (530,000 km$^2$ - larger than Spain and almost as large as France) of north-central and northwestern China, including much or part of seven provinces: Qinghai, Gansu, Ningxia, Inner Mongolia, Shaanxi, Shanxi and Henan, which is shown in Fig. 1 and 2 [2]. The topography of the Loess plateau consists of hilly regions and tablelands extensively cut by rivers and gullies, which can be depicted from Fig. 3 and 4. The result is severe and chronic soil erosion arising from a combination of causes [15].

The climate of the Loess plateau is classified as semi-arid; precipitation varies greatly among areas and years, and tends to fall in high-intensity rainstorms. The region is thus subject to alternating drought and flood, both of which seem to have increased in frequency and severity during recent times. The dense loess soil that gives the plateau its name does not absorb precipitation readily, and most of the rainfall is lost as run-off. Population pressure (currently some 70 million residents, 2011 [16]) has resulted in removal of the native woody vegetation cover, the extent of which has decreased from about 40% cover across the entire region in about 200 A.D. to about 6% in 1980. The result is severe soil depletion over 70% of the region, an annual average topsoil loss of 2500 million tons (much of it carried away by the Yellow River), and a network of large ravines cutting the plateau, including more than 320,000 gullies longer than 1 km, 47,000 ravines between 6 and 10 km in length, and about 2500 longer than 10 km [2].

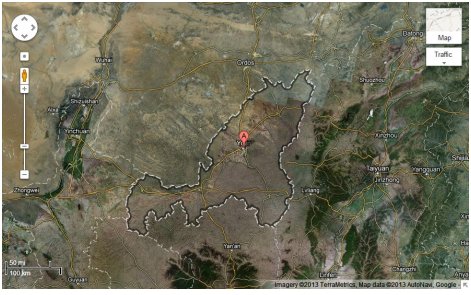

Fig. 5. Graphic view of Yulin region in Loess Plateau from Google Maps. (cited from http://maps.google.com )

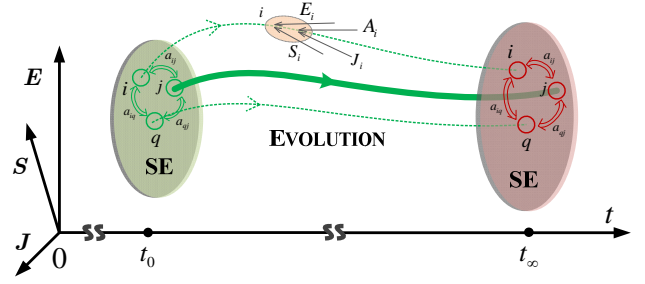

Fig. 6. Diagram of Evolutionary processes of Ecosystem.

The Loess plateau has traditionally supported an agricultural economy dominated by winter wheat. The grain crop is used first to meet the food needs of the farmer and his family. Chaff and some supplementary grain are fed to swine, whereas straw is used for fuel. The farmer's greatest incentive is to maximize production each year, with little regard for long term management. Beyond manuring the wheat fields, little is done to renew the soil. Where poor agriculturalpractices have been combined with soil loss from erosion, the result has been extreme loss of soil fertility and reduction in arability.

In this paper, we will focus our interest in the central part of Loess Plateau, namely *Yulin* region as shown in Fig. 5. According to the data from Chinese Academy of Sciences Water Research Institute, two meters below the huge amount of soil reservoir is the lower bound amount of water to maintain the ecological balance, which however, due to overexploitation, dry soil layer has formed and pose a negative impact on SE planning. Therefore, the optimal water resource management in this region is very important to SE planning of the whole Loess Plateau, which will be focused in Section IV and V.

### III. OVERVIEW OF SE PLANNING OF THE LOESS PLATEAU

#### A. Complexity of evolutionary processes of the Loess Plateau Ecosystem

In exploring SE planning of the Loess Plateau, we face lots of uncertainties which bring diversity, singularity, evolution and complexity. Thus, sustainability evolution can be regarded as nested layers, which not only possess individual-level integration, above-level integration and sectoral-level integration. Moreover, within the same integration level, it also has different subsystems, showing the complex relationships among their levels. Thus, SE is also regarded as a complex giant system. Due to intrinsic uncertainty and randomness of SE, the growth of its complex structure is associated with increment of metabolism of each subsystem [17]. If we treat the continuous evolutionary process of Loess Plateau resources as a combination of all closely related and integrated effective subsystems, in which each subsystem has a certain correlation. Thus, such evolution-based complex system can be formulated as follows.

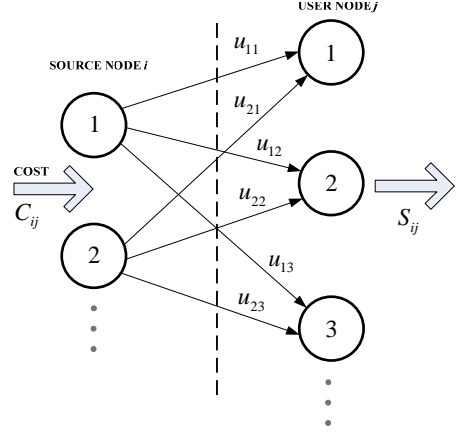

Fig. 7. Diagram of ecological resource planning problem for the Loess Plateau.

$$A(t) = \sum_{i=1}^{n} \left[ (1 + \sum_{j=1}^{n} \sum_{i=1}^{n} a_{ij}) A_i(t) \right] \qquad (1)$$

where $A(t)$ denotes the capability of overall sustainable development of resources in the Loess Plateau; $A_i(t)$ denotes the capability of sustainable development of resources of subsystem $i$; $a_{ij}$ denotes the interaction coefficient between subsystem $i$ and $j$, which reveals the contribution of their interactions to the whole ecosystem; $t$ belongs to the range of $[t_0, t_\infty]$ and $t_\infty$ is our interested time period. $A_i(t)$ can be solved from the following differential equations.

$$\dot{A}(t) = \boldsymbol{f}\left(A_i, \boldsymbol{E}_i, \boldsymbol{S}_i, \boldsymbol{J}_i, t\right) \qquad (2)$$

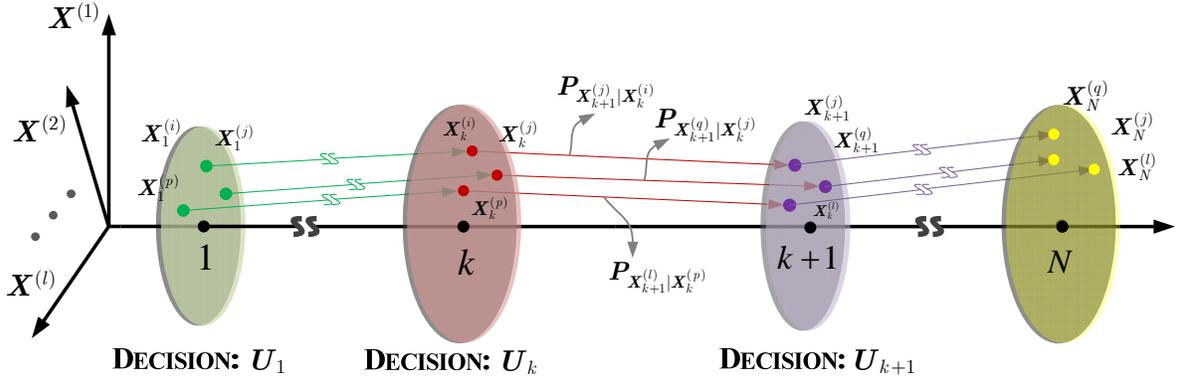

Fig. 8. Illustrative diagram of SE planning based on SDP modeling.

where $E_i$ represents the ecological condition of subsystem $i$; $S_i$ represents the social condition of subsystem $i$; $J_i$ represents the economic condition of subsystem $i$. It should be emphasized that the above vector variables are also time-driven along with the evolution of SE. Thus, an optimization-based formal model can be derived as follows[1].

$$\max_{E_i, S_i, J_i} \Psi = \int_{t_0}^{t_\infty} A(t) dt$$

$$\begin{aligned}
s.t. \quad & \dot{A}_i(t) = f\left(A_i, E_i, S_i, J_i, t\right) \\
& \underline{A}_i^{(k)} \le A_i^{(k)} \le \overline{A}_i^{(k)}, \quad \forall i \\
& \underline{E}_i^{(k)} \le E_i^{(k)} \le \overline{E}_i^{(k)}, \quad \forall i, \forall k \\
& \underline{S}_i^{(l)} \le S_i^{(l)} \le \overline{S}_i^{(l)}, \quad \forall i, \forall l \\
& \underline{J}_i^{(m)} \le J_i^{(m)} \le \overline{J}_i^{(m)}, \quad \forall i, \forall m \\
& t_0 \le t \le t_\infty
\end{aligned} \quad (3)$$

where each decision variable has its upper and lower bound as shown in the above constraints. This model is similar to the idea proposed in [17], which shows that by fully mobilizing the initiative and dynamics of each subsystem, the optimal overall system objective can be achieved through the development and synergies of each subsystem. Fig. 6 illustrates this evolution process.

Thus, the SE planning of Loess Plateau is reliant on evolution rules and focus on the metabolic behaviors such as absorption and regeneration, so as to guide the overall development of ecosystem in the optimal direction. For simplification, we will only focus on the ecological condition $E_i$ of subsystem $i$ and analyze the resource management problem of $E_i$ using SDP in Section IV.

### A. Evolutionary Game Theory

Game theory is a powerful tool to investigate interactions processes among different decision-making agents and their results, which also explore rational behavior under interdependence. There are two basic assumptions: each topic has a clear exogenous variable; each decision from decision-makers is based on their knowledge and expected/prediction of other decision-makers [18].

In recent years, evolutionary game theory has combined the dynamic evolution and game theory, which assumes that there are many participants in a system [19]. Depending on the analytical framework, the number of participants can be countable or uncountable. Each game can be regarded as a random sampling among the set of all participants of game, and the selected participants will do predetermined game (i.e., game elements) so as to get benefits from these game elements. In the framework, we assume that countable participants can't know exactly their current state of interest. Thus, the most favorable strategic will gradually simulate and continue to go, and ultimately achieve a state of equilibrium. It has been shown that in such changeable system, the ratio of people who choose high reward strategy will gradually increase. As one of basic concepts of evolutionary game theory, the *evolutionary stable strategy* (ESS) can be used to simulate a group of ecosystem such as humankinds, animals, plants, whose evolution patterns could be evolutionary stable [19], [20]. In this paper, we will implement this idea explore dynamic decisions of how to cooperate together to achieve sustainable use of resources in SE planning of Loess Plateau.

## IV. SE PLANNING OF LOESS PLATEAU ECOSYSTEM BASED ON SDP

As aforementioned, the ecological resource management such as water, land, plants has become a very important issue in SE planning. The evolution process of Loess Plateau mentioned in Section II has shown that natural and human factors threat the sustainable development of Loess Plateau, especially the shortage of water resource. To reduce negative impact of ecological resource shortage, much effort should be

---

[1] Here, the proposed optimization model is similar to the hybrid system model that I discussed, which can be found in http://arxiv.org/abs/1305.0978. Actually, I am not sure that this is a correct model for ecosystem such as Loess Plateau; in fact, research related to ecosystem is rather difficult as discussed before due to its complex structures and metabolic processes. The purpose of this paper is to explore some feasible applications of decision theory methods into the SE planning.

DRAFT                                                                                                                                        4

Table I  Optimal Result of Water Planning of Yulin Area based on the Proposed Model

| Water Resources | | Agricultural Irrigation | Industrial Usage | Daily Usage |
|---|---|---|---|---|
| Surface Water | Optimal Amount | 2.1307 | 0.334 | 0.175 |
| | Ratio | 44.9% | 64.6% | 29.71% |
| Ground Water | Optimal Amount | 2.6151 | 0.183 | 0.414 |
| | Ratio | 55.1% | 35.4% | 70.29% |

Table II  Actual Result of Water Planning of Yulin Area from Water Data Bulletin [16]

| States of PSS Parameters | | Agricultural Irrigation | Industrial Usage | Daily Usage |
|---|---|---|---|---|
| Surface Water | Actual Amount | 3.5239 | 0.194 | 0.153 |
| | Ratio | 74.5% | 37.3% | 25.9% |
| Ground Water | Actual Amount | 1.2062 | 0.3260 | 0.4372 |
| | Ratio | 25.5% | 62.7% | 74.1% |

focus on the optimal resource management and control. In this Section, we will propose a simplified SDP model based on discrete-time stochastic process for SE planning. Specifically, we will use the proposed model to analyze water resource planning of Loess Plateau.

### A.  Simplified SDP model for SE Planning

First, let's define the total kinds of concerned resource as $m$, which is utilized by $n$ user subsystems (or known as $n$ users). And these users can be regarded as residents, companies, governments, agriculture firms, and etc. For clarity, we give the following definitions of different variables.

*Definition 4.1* (Time Horizon): Define time horizon as $k = 1, 2, ..., N$, which represent the period when each user begin to utilize the resource.

*Remark 4.1*: It is obviously that if $N$ is larger than more accurate decision policies can be obtained. In this paper, we will focus the time horizon in one year, i.e., $N = 12$.

*Definition 4.2* (State Variable): Define state variable at time $k$ as follows

$$\boldsymbol{X}_k = \begin{bmatrix} x_{11}(k) & x_{12}(k) & \cdots & x_{1n}(k) \\ x_{21}(k) & x_{22}(k) & \cdots & x_{2n}(k) \\ \vdots & \vdots & \ddots & \vdots \\ x_{m1}(k) & x_{m2}(k) & \cdots & x_{mn}(k) \end{bmatrix}, \forall k = 1, ..., N \quad (4)$$

where its element $x_{ij}(k)$ denotes as the case of whether the $i$-th resource is used by the $j$-th user. If $x_{ij}(k) = 1$, then the $i$-th resource is assigned to the $j$-th user; if $x_{ij}(k) = 0$, then the $i$-th resource is not assigned to the $j$-th user.

*Remark 4.2*: The state variable $\boldsymbol{X}_k$ belongs to $\mathbb{R}^{m \times n}$ matrix space with its elements of integer 0/1. Mathematically speaking, there are $2^{m \times n}$ kinds of possible selection of $\boldsymbol{X}_k$. However, for practical application, it is obviously that we can't select every element of $\boldsymbol{X}_k$ as zero, which means there is no resource is assigned to any user. In this paper, we assume that each resource will be assigned to arbitrary user; and each use can get at least one kind of resources. Thus, each row and each column of $\boldsymbol{X}_k$ will have at least an integer 1 for any $k = 1, ..., N$. We will see that it is a practical assumption in statistic sense, especially for water resource planning of Loess Plateau.

Moreover, a stationary Markov chain is used to generate the state variable $\boldsymbol{X}_k$, which is assumed to take on a finite number of values

$$\boldsymbol{X}_k \in \{\boldsymbol{X}_k^{(1)}, \boldsymbol{X}_k^{(2)}, ..., \boldsymbol{X}_k^{(i)}, ..., \boldsymbol{X}_k^{(l)}\} \triangleq \tilde{\boldsymbol{X}}_k \quad (5)$$

*Definition 4.3* (Decision Variable): Define decision variable at time $k$ as follows

$$\boldsymbol{U}_k \triangleq \boldsymbol{U}_k(\boldsymbol{X}_k) = \begin{bmatrix} u_{11}(k) & u_{12}(k) & \cdots & u_{1n}(k) \\ u_{21}(k) & u_{22}(k) & \cdots & u_{2n}(k) \\ \vdots & \vdots & \ddots & \vdots \\ u_{m1}(k) & u_{m2}(k) & \cdots & u_{mn}(k) \end{bmatrix} \quad (6)$$

where its element $u_{ij}(k)$ denotes as the case of how much amount of water that the $j$-th user decide to use from the $i$-th resource.

*Remark 4.3*: As aforementioned, $x_{ij}$ can determine whether $j$-th user can utilize the $i$-th resource. Thus, for any $x_{ij} = 0$, we have $u_{ij} = 0$; for any $x_{ij} = 1$, we have $0 < u_{ij} \leq \bar{u}_{ij}$. Then, we can find that $\boldsymbol{U}_k$ is dependent of $\boldsymbol{X}_k$ with a similar matrix structure. We will use this fact in the following discussions. For clearance, let's denote the set of decision variables at $k$-th period as $\tilde{\boldsymbol{U}}_k(\boldsymbol{X}_k)$, where $\boldsymbol{U}_k \subset \tilde{\boldsymbol{U}}_k(\boldsymbol{X}_k)$.

*Definition 4.4* (Transition Probability Matrix): Define transition probability matrix at time $k$ with state $\boldsymbol{X}_k$ and decision $\boldsymbol{U}_k$ as follows

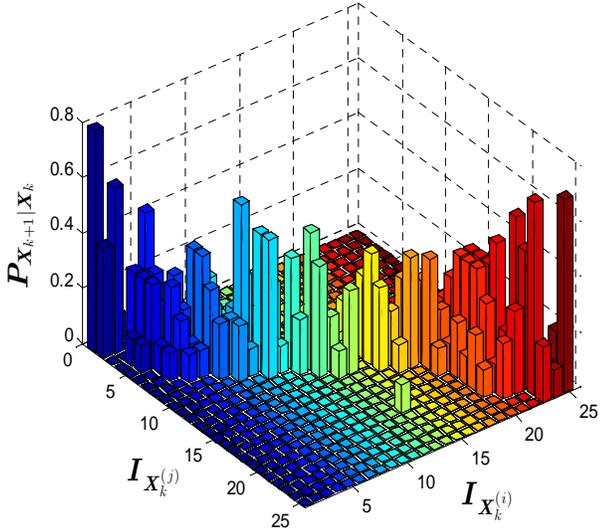

Fig. 9. Histogram of transition probability of state variables at $k=5$.

$$\boldsymbol{P}_{\boldsymbol{X}_{k+1}|\boldsymbol{X}_k} \triangleq \boldsymbol{P}(\boldsymbol{X}_{k+1} \mid \boldsymbol{X}_k, \boldsymbol{U}_k)$$
$$= \begin{bmatrix} p_{11}(k) & p_{12}(k) & \cdots & p_{1l}(k) \\ p_{21}(k) & p_{22}(k) & \cdots & p_{2l}(k) \\ \vdots & \vdots & \ddots & \vdots \\ u_{l1}(k) & u_{l2}(k) & \cdots & p_{ll}(k) \end{bmatrix} \quad (7)$$

where $p_{ij}(k)$ denotes transition probability from $\boldsymbol{X}_k^{(i)}$ to $\boldsymbol{X}_k^{(j)}$ at time $k$, and $\sum_{j=1}^{l} p_{ij}(k) = 1$.

*Remark 4.4*: As indicated in Definition 4.3, we can simplify $\boldsymbol{P}(\boldsymbol{X}_{k+1} \mid \boldsymbol{X}_k, \boldsymbol{U}_k)$ as $\boldsymbol{P}(\boldsymbol{X}_{k+1} \mid \boldsymbol{X}_k)$. However, the transition probability matrix $\boldsymbol{P}_{\boldsymbol{X}_{k+1}|\boldsymbol{X}_k}$ is a statistic matrix which can't be easily derived from analytical modeling. In this paper, we will use the statistic data of water resource bulletin [16] to determine $\boldsymbol{P}_{\boldsymbol{X}_{k+1}|\boldsymbol{X}_k}$. By using the maximum likehood estimator, $\boldsymbol{P}_{\boldsymbol{X}_{k+1}|\boldsymbol{X}_k}$ could be estimated as the observation data as follows:

$$\hat{p}_{ij}(k) = \frac{\tilde{N}_{ij}}{\tilde{N}_i}, \quad \forall k = 1,...,N \quad (8)$$

where $\tilde{N}_{ij}$ is the number of occurrences of the transition from $\boldsymbol{X}_k^{(i)}$ to $\boldsymbol{X}_k^{(j)}$ at time $k$, and $\sum_{j=1}^{l} \tilde{N}_{ij} = N_i$ is the total number of times that $\boldsymbol{X}_k^{(i)}$ has occurred at time $k$. However, it is possible a training set from water resource bulletin may not be rich enough to cover the whole state space, and, in some cases, $N_i$ may be zero form some $\boldsymbol{X}_k^{(i)}$. Therefore, a smoothing technique was used to the estimated parameters [21], [22].

*Definition 4.4* (Reward Function): Define reward function at time interval $[k, k+1]$ as follows

$$\boldsymbol{V}_k = \begin{bmatrix} r_{11}(k) & r_{12}(k) & \cdots & r_{1n}(k) \\ r_{21}(k) & r_{22}(k) & \cdots & r_{2n}(k) \\ \vdots & \vdots & \ddots & \vdots \\ r_{m1}(k) & r_{m2}(k) & \cdots & r_{mn}(k) \end{bmatrix} \quad (9)$$

where the element $r_{ij}(k)$ can be express as

$$r_{ij}(k) = \begin{cases} S_{ij} - C_{ij} \cdot u_{ij}(k), & x_{ij}(k) \neq 0 \\ 0, & x_{ij}(k) = 0 \end{cases} \quad (10)$$

where $S_{ij}$ denote the reward of the *j*-th user that utilized the *i*-th resource; $C_{ij}$ denote the cost of the *j*-th user that utilized per-unit amount of the *i*-th resource. We can assume that the reward of using each resource for the same user and the cost for using the same resource for different users remain unchanged with respect to time $k$, respectively.

Therefore, we can simplify (10) as follows

$$r_{ij}(k) = \begin{cases} S_{\cdot|j} - C_{i|\cdot} \cdot u_{ij}(k), & x_{ij}(k) \neq 0 \\ 0, & x_{ij}(k) = 0 \end{cases} \quad (11)$$

where $S_{\cdot|j}$ and $C_{i|\cdot}$ can be also found from the statistic data of water resource bulletin [16].

Fig. 7 and 8 give intuitive descriptions of the above SE planning problems along with these definitions and remarks. Finally, to solve this planning problem, we can use DP algorithm as follows [23]

$$\begin{aligned} &J_k(\boldsymbol{X}_k) \\ &= \max_{\boldsymbol{U}_k \in \tilde{\boldsymbol{U}}_k} \left\{ \sum_{i=1}^{m} \sum_{j}^{n} r_{ij}(k) + \sum_{\forall \boldsymbol{X}_k \in \tilde{\boldsymbol{X}}_k} \boldsymbol{P}_{\boldsymbol{X}_{k+1}|\boldsymbol{X}_k}(i,j) g_{k+1}(\boldsymbol{X}_{k+1}) \right\} \\ &= \max_{\boldsymbol{U}_k \in \tilde{\boldsymbol{U}}_k} \left\{ \sum_{i=1}^{m} \sum_{j}^{n} r_{ij}(k) + \sum_{\forall \boldsymbol{X}_k \in \tilde{\boldsymbol{X}}_k} p_{ij}(k) J_{k+1}(\boldsymbol{X}_{k+1}) \right\} \end{aligned} \quad (12)$$

where $J_k(\boldsymbol{X}_k)$ denote the total reward within time interval $[k, N]$, and the boundary condition $J_N(\boldsymbol{X}_N)$ and $\boldsymbol{X}_N$ are known according to statistic data.

### B. Analysis of Water Resource Planning based on the Proposed Model

For the water resource planning in Loess Plateau, we would like to address the following facts:
1) The water resources can be classified as two parts: surface water and ground water, i.e., $i = 2$;
2) The users subsystems can be classified as three parts: agricultural firms, industrial usage and daily usage, i.e., $j = 2$
3) As indicated in 2011 Water Data Bulletin of Shannxi Province (2011 year) [16], $S_{\cdot|j}$ and $C_{i|\cdot}$ in (11) can be selected as

$$S_{\cdot|j} = \begin{bmatrix} \cdots \\ 3 & 4 & 1 \\ \cdots \end{bmatrix}, \qquad C_{i|\cdot} = \begin{bmatrix} \vdots & 3 & \vdots \\ & 5 & \end{bmatrix} \quad (13)$$

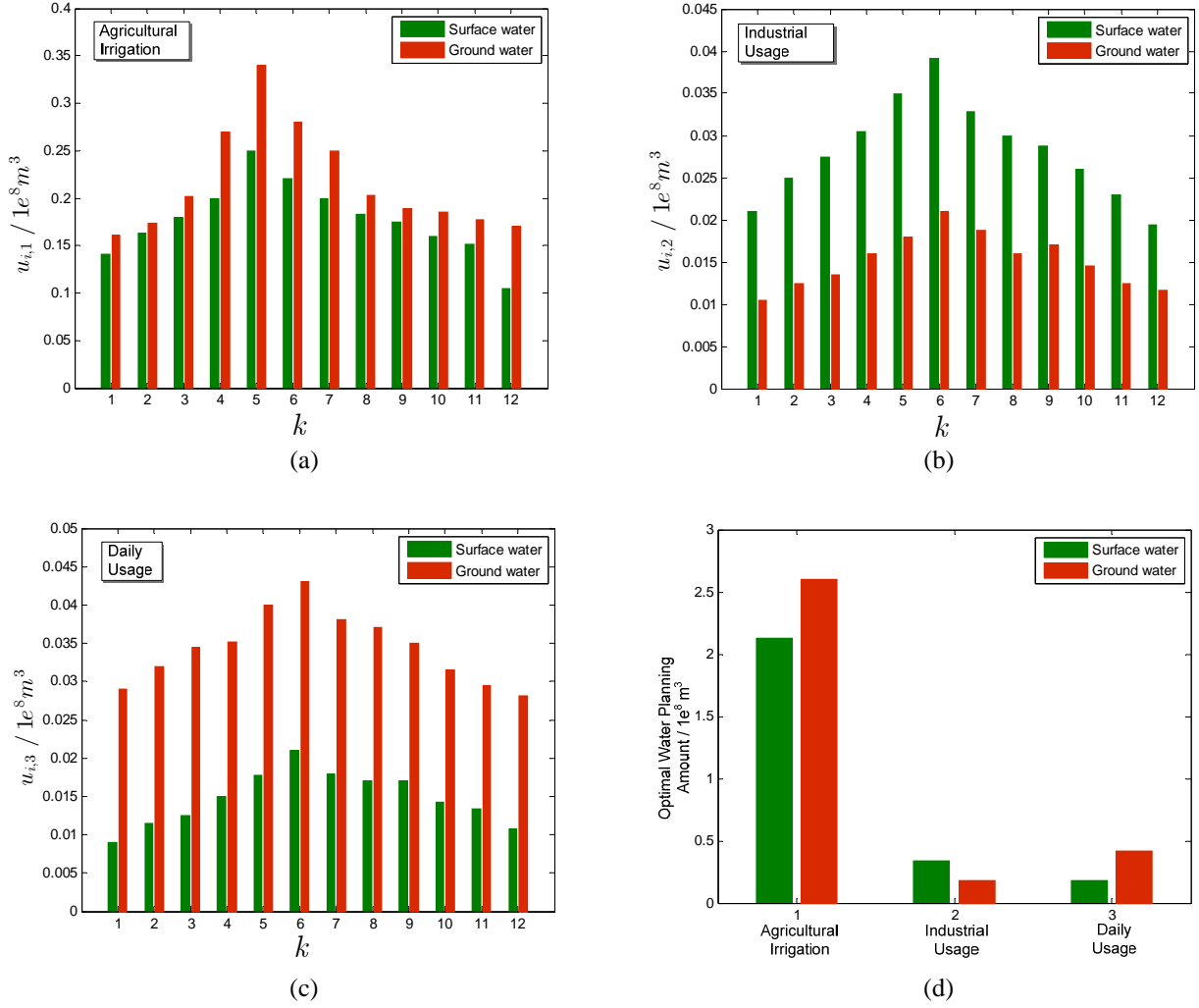

Fig. 10. Optimal results of water planning of Yulin area. (a), (b), (c) Per-month amount of agricultural irrigation, industrial usage, daily usage, respectively; (d) Whole year optimal amount of water planning.

4) As indicated in Water Data Bulletin of Shannxi Province [16], the water resource for Loess Plateau are main restricted into the above two parts. In our water planning, we expect to fully and sustainable use both of water resources. As mentioned, the case for each row and each column of $\boldsymbol{X}_k$ have at least an integer 1 for any $k=1,...,12$ should be a reasonable assumption. Based on it, we list all the 25 possible cases in Appendix A, where the index set $\boldsymbol{I}_{\boldsymbol{X}_k^{(i)}} = \{1,2,...,i,...,25\}$ maps each state matrix variable $\boldsymbol{X}_k^{(i)}$ into the corresponding integer $i$. Thus, the transition probability matrix $\boldsymbol{P}_{\boldsymbol{X}_{k+1}|\boldsymbol{X}_k}$ in (7) can be rewritten as

$$\boldsymbol{P}_{\boldsymbol{X}_{k+1}|\boldsymbol{X}_k} = \begin{bmatrix} p_{11}(k) & p_{12}(k) & \cdots & p_{1,25}(k) \\ p_{21}(k) & p_{22}(k) & \cdots & p_{2,25}(k) \\ \vdots & \vdots & \ddots & \vdots \\ p_{25,1}(k) & p_{25,2}(k) & \cdots & p_{25,25}(k) \end{bmatrix}, \quad (14)$$

$$\forall k = 1,...,12$$

Thus, by using (8), we can obtain (14) at each time $k$. Fig. 9 depicts the transition probability matrix $\boldsymbol{P}_{\boldsymbol{X}_{k+1}|\boldsymbol{X}_k}$ as $k=5$. Moreover, combining (11), (12), (13), (14), we can solving the SDP model and get the optimal water planning policies for each time $k=1,...,12$. Fig. 10 Shows the optimal results of water planning of Yulin area, where the whole year amount of water depicted in Fig. 10 (d) corresponds with data of Table I. Compared with the actual amount of water planning shown in Fig. 11 (corresponding with data of Table II), we find that the actual water resource planning is not reasonable. Several interim conclusions on the water planning strategies could be drawn based on the above numerical results.

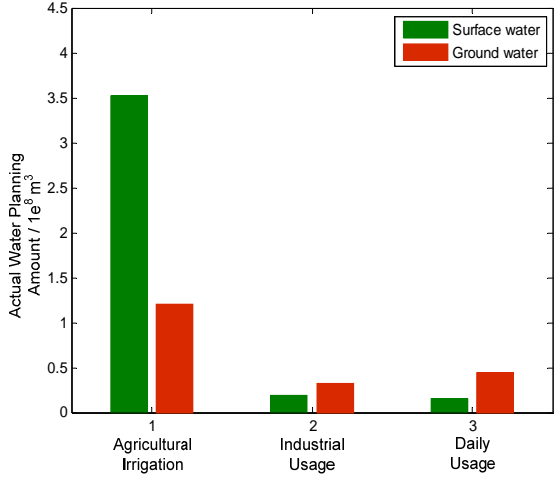

Fig. 11. Actual amount of water planning of Yulin area from 2011 Water Data Bulletin [16].

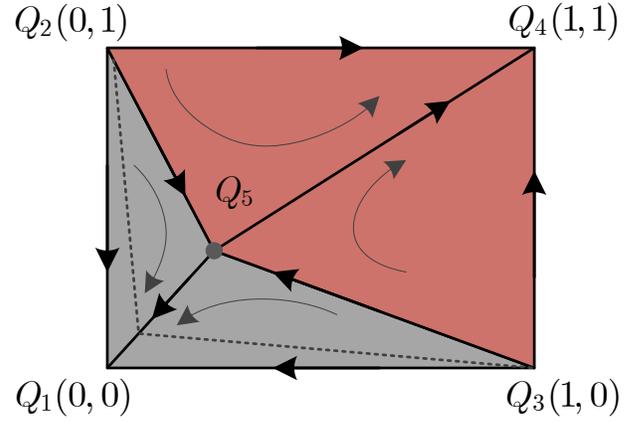

Fig. 12. Dynamic evolutionary process of equilibrium in the SE planning.

1) At the current stage, agricultural irrigation should be focus on the utilization of ground water, and enlarge its the amount ratio;
2) The industrial and daily usage of water should focus on surface water.

Some other detailed suggestions such as water facility constructions could be naturally obtained based on the above two points, which are not the focus of interests herein.

## V. Dynamic Game Analysis of Water Resource Planning of Loess Plateau

Section IV focuses on the analysis of overall water resource planning problem among multi-users, which however, ignores impacts of dynamic relationship among these users, namely evolutionary game relations, such as the influence of polluted water from industrial to agricultural and daily usage water. In this section, we will apply the evolutionary game theory as a supplementation of the proposed SDP model. It is worth emphasizing that though our analysis is focus on water planning, the fundamental idea is also applicable to other SE planning of Loess Plateau.

First, assume that there contain $n$ number of social members, known as $n$ participants in the water resource planning. Suppose that there are two strategies in the decision games, namely, $C$ (sustainable usage), $D$ (unsustainable usage). At the time $t=0$, each participant randomly choose a strategy. Then for time $t=1$, only one participant can change his/her previous decision [19]. Also assume that we randomly pick two participants to do game playing in group $A$ and $B$, where the two participants are independent with each other if choosing properly. In this case, each payoff equals to 1, and $u$, $v$ ($u>1, v>1$) denote as the payoff of $A$ and $B$, under cooperation case, respectively.

Denote $p$ as the ratio of participant who choose strategy C among group A; $q$ as the ratio of participant who choose strategy D among group B. Thus, $(p,q)$ can represent the evolutionary dynamics of the system. As indicated in [19], if the adaptability of a strategy is higher than the average adaptability of a group, then such strategy can be developed. The dynamics of the ratio $p$ and $q$ can be written as

$$dp/dt = p(1-p)(uq-1) \qquad (15)$$

$$dq/dt = q(1-q)(vp-1) \qquad (16)$$

where equilibrium $p=0;1$ or $q=1/u$ can lead to stable ratio of participant who choose strategy C in group A. Likewise, equilibrium $q=0;1$ or $q=1/v$ van lead to stable ratio of participant who choose strategy D in group B. Though calculating the determinants and traces of Jacobia matrix of (15), (16) in the neighborhood of local equilibrium, we can analyze the pattern of local equilibrium, which is shown in Table III.

Fig. 12 depicts the dynamic evolutionary process of the above five equilibrium in SE planning. Here, only $Q_1 \triangleq (0,0)$ and $Q_4 \triangleq (1,1)$ are stable, which belongs to ESS as indicated in [1]. According to (15), (16), we can see that they also correspond to the sustainable mode and over-exploitation one. The polygonal lines formed by unstable points $Q_2 \triangleq (0,1)$, $Q_3 \triangleq (1,0)$, $Q_5 \triangleq (1/v, 1/u)$ can be regarded as boundaries of different states.

It is clearly that all the points at the lower-left part of polygonal line $Q_2Q_5Q_3$ will eventually converge to $Q_1$, namely over-exploitation mode; other points at the upper-right part will converge to $Q_4$, namely sustainable mode. Therefore, if chosen larger $u$, $v$, then $Q_5$ will converge closer to $Q_1$, which decrease the lower-left part area. Thus, the probability that the system states will have larger opportunities to converge to $Q_4$ rather than converging to $Q_1$, which is more prone to be sustainable mode. This leads to the idea of increasing payoff so as to achieve a better water/SE planning of Loess Plateau.

Table III  Patterns of equilibrium

| Equilibrium $(p, q)$ | det $(J)$ | | tr$(J)$ | | Results |
|---|---|---|---|---|---|
| (0, 0) | 1 | + | -2 | − | **ESS** |
| (0, 1) | $u$-1 | + | $u$ | + | Unstable |
| (1, 0) | $v$-1 | + | $v$ | + | Untable |
| (1, 1) | $(u-1)(v-1)$ | + | $2-u-v$ | − | **ESS** |
| $(1/v, 1/u)$ | $-(1-1/v)(1-1/u)$ | − | 0 | 0 | Saddle point |

## VI. CONCLUSION

In this paper, the SE planning of the Loess Plateau area in Northwest China has been analyzed based on SDP model and evolutionary game theory. The concept of SE planning is introduced with specifications in ecological resource planning. Detailed mathematical discrete SDP model is formed with several reasonable assumptions, and the transition probability matrix is calculated in a statistic sense. A typical Yulin area in the central part of the Loess plateau is selected to exemplify the validity of proposed SDP model. By implementing a backward DP recursion, the optimal water resource planning strategies are obtained. Compared with actual resource planning condition, several suggestions are generated. Besides quantifying the optimal SE planning problems, a supplementary discussion of impacts between different users is demonstrated in the sense of evolutionary game playing. Although the approach is applied to the water resource planning of Loess Plateau as an example, the methodology of using SDP and game theory is applicable to other ecosystems.

## ACKNOWLEDGEMENT

I would like to thank Y. X. Li (Boston University) for drawing the attention to [21] and [22]. I also wish to thank Grace Dai for her kind supports at different stages of preparing this paper. My deepest gratitude goes to Prof. Michael Caramanis, for his patience and helpful suggestions during my first-academic year of graduate studies at Boston University.

## APPENDIX

Explicit forms of all possible state variables $X_k$ can be listed as follows

$$X_k^{(1)} = \begin{bmatrix} 0 & 0 & 1 \\ 1 & 1 & 0 \end{bmatrix}, X_k^{(2)} = \begin{bmatrix} 0 & 0 & 1 \\ 1 & 1 & 1 \end{bmatrix},$$
$$X_k^{(3)} = \begin{bmatrix} 0 & 1 & 0 \\ 1 & 0 & 1 \end{bmatrix} \quad (17)$$

$$X_k^{(4)} = \begin{bmatrix} 0 & 1 & 1 \\ 1 & 0 & 0 \end{bmatrix}, X_k^{(5)} = \begin{bmatrix} 0 & 1 & 1 \\ 1 & 0 & 1 \end{bmatrix},$$
$$X_k^{(6)} = \begin{bmatrix} 0 & 1 & 0 \\ 1 & 1 & 1 \end{bmatrix} \quad (18)$$

$$X_k^{(7)} = \begin{bmatrix} 0 & 1 & 1 \\ 1 & 1 & 0 \end{bmatrix}, X_k^{(8)} = \begin{bmatrix} 0 & 1 & 1 \\ 1 & 1 & 1 \end{bmatrix},$$
$$X_k^{(9)} = \begin{bmatrix} 1 & 0 & 0 \\ 0 & 1 & 1 \end{bmatrix} \quad (19)$$

$$X_k^{(10)} = \begin{bmatrix} 1 & 0 & 1 \\ 0 & 1 & 0 \end{bmatrix}, X_k^{(11)} = \begin{bmatrix} 1 & 0 & 1 \\ 0 & 1 & 1 \end{bmatrix},$$
$$X_k^{(12)} = \begin{bmatrix} 1 & 1 & 0 \\ 0 & 0 & 1 \end{bmatrix} \quad (20)$$

$$X_k^{(13)} = \begin{bmatrix} 1 & 1 & 1 \\ 0 & 0 & 1 \end{bmatrix}, X_k^{(14)} = \begin{bmatrix} 1 & 1 & 0 \\ 0 & 1 & 1 \end{bmatrix},$$
$$X_k^{(15)} = \begin{bmatrix} 1 & 1 & 1 \\ 0 & 1 & 0 \end{bmatrix} \quad (21)$$

$$X_k^{(16)} = \begin{bmatrix} 1 & 1 & 1 \\ 0 & 1 & 1 \end{bmatrix}, X_k^{(17)} = \begin{bmatrix} 1 & 0 & 0 \\ 1 & 1 & 1 \end{bmatrix},$$
$$X_k^{(18)} = \begin{bmatrix} 1 & 0 & 1 \\ 1 & 1 & 0 \end{bmatrix} \quad (22)$$

$$X_k^{(19)} = \begin{bmatrix} 1 & 0 & 1 \\ 1 & 1 & 1 \end{bmatrix}, X_k^{(20)} = \begin{bmatrix} 1 & 1 & 0 \\ 1 & 0 & 1 \end{bmatrix},$$
$$X_k^{(21)} = \begin{bmatrix} 1 & 1 & 1 \\ 1 & 0 & 0 \end{bmatrix} \quad (23)$$

$$X_k^{(22)} = \begin{bmatrix} 1 & 1 & 1 \\ 1 & 0 & 1 \end{bmatrix}, X_k^{(23)} = \begin{bmatrix} 1 & 1 & 0 \\ 1 & 1 & 1 \end{bmatrix},$$
$$X_k^{(24)} = \begin{bmatrix} 1 & 1 & 1 \\ 1 & 1 & 0 \end{bmatrix}, X_k^{(25)} = \begin{bmatrix} 1 & 1 & 1 \\ 1 & 1 & 1 \end{bmatrix} \quad (24)$$